\documentclass{article}
\usepackage{amssymb,amsmath,theorem,euscript}

 \usepackage{hyperref}

\usepackage[cp1251]{inputenc}
\usepackage[english,russian]{babel}

\usepackage{graphicx}

\begin{document}
	
	\def\sm{\smallskip}

\begin{center}
	\bf\Large
Nikolay Konstantinov and the Konstantinov System 

\bigskip

	\sc \large
	Yury A. Neretin
\end{center}	
 
{\small
The 'Konstantinov System' was a non-standard educational  institution  created by the great
mathematical educator Nikolay Konstantinov (1932-2021), this 'System'  worked (mainly in Moscow) in 1960-80s. 
We discuss some sides of  technologies of the 'System': methods of teaching, organization, approaches to personal growth. In Russian.

\medskip

 Published in Matematicheskoe obrazovanie [Mathematical education],   2021, 3 (99) [О Константиновской системе и беломорских стройотрядах]%
\footnote{Here we add a photo at the end of the paper, references \cite{K1}--\cite{DK2} to technological publications
	of N.N.Konstantinov, and   links to MR Author ID's 
	in the caption of Fig.2 	(without an access to MathSciNet they contain a little information), also we added links to some authors in the bibliography.}. 
 
 \bigskip
 
 The paper originally was written for biologists for the site of the White Sea Biological Station of Moscow State University%
\footnote{Web feed of WSBS MSU, 07.07.2021} --- Konstantinov also had an  influence to school biological education in Moscow, on the other hand brigades of co-workers and pupils  of the 'System' made a substantial contribution to the building of this Arctic marine Station.}
 
\bigskip

 \begin{center}
 
 \Large\bf О  Константиновской системе\\ и Беломорских стройотрядах	
 
 \bigskip
 
 \large\sc Ю.А.Неретин
 	
 \end{center}

\bigskip

{\small В статье рассказывается о <<Константиновской системе>> -  необычной общественной образовательной структуре,
	созданной Николаем Николаевичем Константиновым (1924-2021) и работавшей в математическом образовании 60-80~гг. XX~в., а также о беломорских математических стройотрядах 1969-1987~гг., представлявших из себя своего рода ответвление этой системы.}

\newpage
 
3 июля 2021 г. на 90-м году жизни от ковида скончался Николай Николаевич Константинов,  сыгравший видную роль в истории Беломорской биологической станции МГУ%
\footnote{Эта статья первоначально была написана для сайта Беломорской биологической станции МГУ им. Н.А.Перцова (ББС МГУ), сайт {\tt http://wsbs-msu.ru/}, лента за 07.07.2021.  
	\newline
	ББС - морская биостанция, расположена на северной границе  Карелии на Полярном круге на мысу Киндо. На картах начиная с 60х годов отмечается как <<поселок Приморский>>. Ближайшая ж.д. станция - Пояконда, 15 км по прямой. Связь биостанции с миром - по морю малыми судами и моторками (в начале 50х еще и гребными лодками), зимой - по зимнику (лошадью в 50х, потом также на тракторе или грузовике, сейчас - буранами).
	\newline
	Основана в 1938г. проф. Л.~А.~Зенкевичем (был вкопан заявочный столб). В следующем году там уже проходила практику группа студентов и работала водолазная станция (с оборудованием, выданным ЭПРОНом). В Войну биостанция оказалась в прифронтовой полосе, там находился авиационный наблюдательный пост.
\newline
	Летом 1951г. новым директором стал  
	Николай Андреевич Перцов (1924-1987), только что окончивший биофак МГУ. Фронтовик, участник боев под Москвой. Туберкулезник. В момент его приезда  станция представляла из себя избу местного жителя (он же сторож) и несколько дощатых построек сарайного типа в тяжелом состоянии. К середине 70х Перцов превратил ее в знаменитую научно-учебную базу (которая, в частности, принимала ежегодно на практику сотни студентов), см. \cite{Per}, \cite{Kal}. В определенном смысле \cite{X}, биостанция начала представлять из себя и культурный феномен. 
	\newline
	Поселок Приморский строил сам себя (при сильном недостатке финансирования) и  работал как небольшая автономная строительная организация. Перцов совмещал в себе роль директора научно-образовательной организации 
	с фактической должностью хозяйственника и руководителя строительных работ. 
	С самого начала он подключал к работам студентов (местность там красивая, а Перцов, будучи ярким деятелем, умел привлекать к себе людей). Вскоре это превратилось в стройотряды, которые сначала состояли из студентов-биологов (по-видимому, больше биологинь), а потом, вслед за Константиновым, там появились математики и физики.}. 

Я довольно много сотрудничал с ним с 1973-74~гг. по 1980~г., и отчасти вплоть  до 1991-93~гг.,	и хотел бы рассказать
об этом выдающемся человеке и его деятельности в той степени, в которой я мог ее наблюдать (с некоторыми экстраполяциями, основанными на его рассказах и рассказах других людей, сотрудничавших с ним). Начну с событий, которые я не наблюдал и наблюдать не мог. 	

\begin{center}
	\bf
\dots\dots\dots
\end{center}
			
На Беломорскую биостанцию (ББС) Константинов (<<Конст>>, <<НикНик>>) приехал в стройотряд  в июле 1967~г.  
C ним были два его
бывших ученика по математическому классу,  в тот момент уже закончившие первый курс мехмата МГУ. Он рассказывал какую-то запутанную историю о том, как это получилось \cite{Kal}, вероятно его приезд был связан с проблемами со здоровьем у одного из этих учеников.

\begin{figure}	
	\includegraphics[width=157mm]{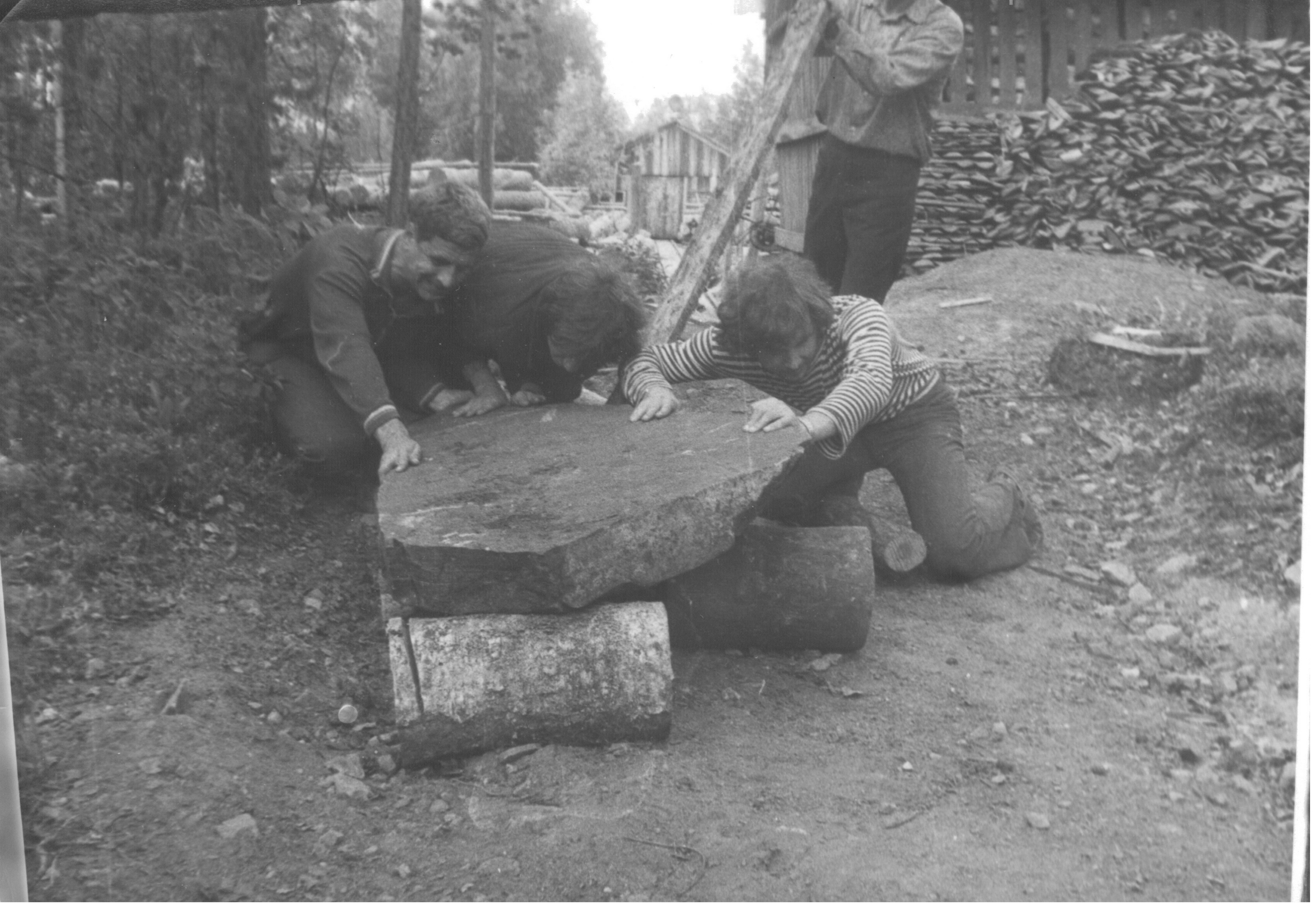}
	\caption{Поселок Приморский (ББС МГУ), 1979~г. Установка <<журнального столика>>  из гнейса  около здания <<Огарки>>.  Затеявший это мероприятие (такие мероприятия назывались <<маразмами>> или <<мехматскими маразмами>>) Константинов - слева. Вдали  - здание пилорамы (с торца напоминает сарай). Остальные объекты на втором плане сейчас не существуют. Слева - штабеля бревен, спиленных и вывезенных из леса  зимой рабочими. Стройотряд вскорости увезет бревна 
		по  рельсам на пилораму, где   превратит их в брусья, плахи, доски разного сечения, бруски,
		рейки, опилки (которые будут собраны в специальном бункере), дрова (они отправятся сохнуть в сарай, который на фотографии справа), а также стружки (которые пойдут в отход).
		Журнальный столик долго стоял на ножках (коротких толстых чурках), в настоящее время не существует.
		\newline	
		Стоит отметить, что на снимке все в порядке с техникой безопасности - в присутствии Конста не могло быть иначе. Камень полностью опирается на катки; силы, прилагаемые к столику участниками транспортировки, направлены по горизонтали в сторону фотографа (с небольшой составляющей вверх, идущей от хиловатой ваги). Кисти рук (кстати, единица масштаба)
		находятся сверху, их не прищемит. Соскользнуть назад или вбок на людей этот столик весом в несколько сот килограмм в запечатленный момент не может. Соскальзывание вперед возможно. Если это случится, то придется ломами и вагами снова водружать столик на катки, но сейчас непосредственной опасности такое соскальзывание ни для кого не несет.	}	
\end{figure}

\begin{quotation}
	\it 
	Еще в Москве мы пообщались с Перцовым, организовывали погрузку вагона\footnote{Каждую весну из Москвы в Пояконду отправлялся вагон с разнообразными грузами.}. Он произвел на меня очень хорошее
	впечатление. Я был, конечно,  старше многих стройотрядовцев\footnote{Константинову было 35 лет.}. Тогда Бурковский%
	\footnote{И.В.Бурковский, гидробиолог, 1942г. рождения, в 1966г. окончил биофак МГУ.} был руководителем стройотряда, и отряд почти целиком состоял из девочек  - учениц вечернего отделения биофака.
	Вот мы немножко и разбавили их.
\end{quotation}

Константинов был человеком рукастым. Он рассказывал  много историй о своих поездках
по стране, был он в частности в стройотряде на Сахалине, который занимался лесоповалом. Я, например, помню, как он нам объяснял, как надо валить тонкие 15-сантиметровые сосны тремя ударами (лесорубского) топора.

Приехавшие на ББС математики  пришлись ко двору и как работники,  и как люди, способные понимать, как организовывать работы. Константинов рассказывал \cite{Kal}, что один из этих студентов посмотрел, как стройотрядовцы перегружают кирпич (тогда строили  Аквариалку -- трехэтажный кирпичный лабораторный корпус), передавая его по цепочке из рук в руки, и объяснил, что лучше перемещать кирпич по цепочке,
{\it кидая} его в руки соседа. Конечно, это была известная технология, но на биостанции ее то ли не знали, то ли боялись применять - она требовала определенной организации и жесткого исполнения определенных мер  техники безопасности.   Сам
Константинов вскоре оказался руководителем разгрузочно-погрузочных работ в Пояконде. 

На следующий, 1968~г. в опубликованных \cite{Kal} списках стройотрядов видно несколько математиков. Константинов там
не значится, но  некоторые рассказчики утверждали, что он приезжал. В 1969~г. Перцову были нужны  рабочие руки  для строительства   линии электропередач Пояконда--ББС. Для тяжелых работ в лесу и на болоте желательны были люди мужеска пола, а их на биофаке не хватало. Тогда и приехала первая большая группа юных математиков (многие из них впоследствии стали известными людьми). С этого пошли  математические стройотряды 
1969-1987~гг. на ББС,  немало там построившие и многое делавшие для жизнеобеспечения биостанции \cite{Per}, \cite{Kal}.

\begin{center}
	\bf
	\dots\dots\dots
\end{center}

 Здесь надо прерваться и напомнить, что в другой своей ипостаси Константинов   был, быть может, самым ярким деятелем  нашего школьного математического образования за последние сто лет. Попытаюсь рассказать о нем  без восклицательных знаков и эпитетов.

\sm 

Он был%
\footnote{См. о нем \cite{Mat}.} одним из организаторов первых математических школ 
(начало 60х) и уже в 60х годах создал так называвшуюся <<Константиновскую систему>> - общественную структуру в Москве,
которая контролировала несколько математических школ (в разное время 444, 7, 57, 91 и 179, столь же знаменитые  2-ая школа и Колмогоровский интернат в эту систему не входили), сеть математических кружков, математические  олимпиады, летние математические лагеря для школьников и студентов (список этот заведомо не полон). Как известно, бывает немало такого, что существует на бумаге, не существуя в жизни. Константиновская система существовала только в жизни, на бумаге же ее не было.

Люди, начинающие учить сильно продвинутых школьников математике, сталкиваются с неочевидной дилеммой, чему учить -  <<углубленной элементарной математике>> или <<высшей математике>>, которой учат в институтах/университетах. Первый путь ведет к развитию не очень естественных приемов решения задач и искусственному раздуванию их сложности (думаю, что многие биологи не лишены такого школьного опыта или наблюдали его у своих детей). Зачем идти вторым путем не очень ясно (особенно, если школьники настроены поступать в вузы, в которых это все равно будут проходить на более высоком уровне). Варианты этой проблемы должны возникать и при обучении продвинутых школьников другим наукам. 

Константинов нашел нетривиальное решение для математики в виде <<бега на месте>>. В школах <<Системы>> был введен дополнительный предмет <<мат.анализ>>.  На нем решались интересные элементарные задачи, а
также (это занимало основное время) проходилась за три года тогдашняя программа мат.анализа мехмата МГУ ... за один месяц первого семестра (и проходилась, скорее, не полностью). Школьникам предлагались <<листочки>>, отпечатанные на папиросной бумаге (на пишущей машинке под копирку), с длинными сериями задач, которые ученик должен был решить сам.
Решение принималось у каждого по отдельности в письменном виде, и его надо было объяснить в ходе пренаиподробнейшего разговора с потоком неожиданных вопросов со стороны принимающего. Что касается классной доски, то она использовалась разве что для написания объявлений. Учитель публично не объяснял ничего (sic!). Все делалось мягко, уровень педагогического насилия был почти нулевой (если человек <<увязал>>, ему старались помочь боковыми вопросами) ... В итоге за 2-3 года почти ничего не проходилось, но ... учащийся приобретал навыки самостоятельного мышления, точных логических рассуждений и интеллектуальной инициативы... На это и делалась ставка.
Были определенные минусы - у школьника могло возникать  впечатление,  что теорию пределов он чуть ли не придумал сам (то есть он столь же велик, как такие классики как Евдокс, Больцано или Коши). Конечно,
это было иллюзорным - в действительности его вели малыми шагами по тщательно проложенным и ухоженным тропинкам. Были и иные минусы, но плюсы сильно перевешивали.

На эту же цель - учить думать, а не накачивать информацией - были настроены и кружки для школьников (это легко сказать, но не столь уж просто сделать, нужны были подходящие задачи, нужна была технология проведения занятий, нужны были люди, готовые играть в такие игры, была нужна обеспечивающая приход школьников реклама...).  

Математические олимпиады времен <<Системы>> следовали традициям математических олимпиад (они
появились еще до Войны) и были мало похожи на большинство современных олимпиад, так что одно из этих слов
- старое или современное - следовало бы брать в кавычки. Разумеется, это было соревнование, но реальных <<пряников>> на олимпиадах почти не раздавалось  (это  была установка Министерства). Олимпиады по математике (как и по прочим предметам) были просветительским мероприятиями - показать молодежи, как красивы науки и выявить способных людей, которые могут заинтересоваться ими (быть может, этот дух остался в ШБО - Школьной биологической олимпиаде биофака МГУ, которая, впрочем, и не дает никаких пряников). 

Кстати, в олимпиадах, курировавшихся лично Константиновым или <<Системой>>,  была жесткая установка - задача бывает либо решенной, либо не решенной. Всё прочее - от лукавого. Никаких <<подводных камней>> и учета недочетов. Задача <<в общем решена>> - почти то же, что решена. <<В общем не решена>> - почти то же, что не решена (ошибки в работах, конечно, тщательно выявлялись и указывались, но реально разница между <<почти решена>> и <<вполне решена>> могла иметь значение лишь в единичных экзотических случаях). Для  тех, кто в последние десятилетия подвергался оцениванию или сам оценивал, это повод задуматься о разных системах оценивания письменных работ и влиянии этого оценивания на обучаемого и обучение.  Кстати, проверка была делом непростым - в решениях требовалось развернутое логическое обоснование, к чему школьники, разумеется, были не приучены, а задачей проверяющего было доброжелательно разобраться в этой написанной некалиграфическим почерком абракадабре...

\medskip 

Читатель может спросить - а кто же этим занимался? А занималась большая веселая тусовка, состоявшая, в основном, из студентов и аспирантов\footnote{А в ведении кружков было обычным участие старшеклассников.}. В основном из тех, кто раньше учился в <<Системе>>. На школьном уроке
одновременно присутствовали несколько преподавателей, и все они одновременно беседовали со школьниками,
переходя от одного к другому. Мертвящей тишины в классе - ой - не было. Разумеется, удельные трудовые затраты были
очень велики, но участникам эта живая игра была интересна, да и они сами чему-то обучались. А когда (в норме - не если, а когда) становилось скучно, или когда людей одолевали другие проблемы - они просто уходили из игры, это считалось
само собой разумеющимся. Тем временем подрастали новые молодые люди, исходившие из того, что с другими надо делать то, что делали раньше с тобой. Как ни странно это может выглядеть издалека, непосредственная работа держалась на непрофессионалах и непрофессионализме. Предполагалось, что они все делают сами и сами решают, что делать. 

Во время, которое я наблюдал, школьные учителя были вне <<Системы>> и находились с ней в симбиозе. С ними не могли не контактировать, но у них были свои предметы, у <<Системы>> - только матан. Среди учителей были
выдающиеся люди (как В.М.Сапожников или В.В.Бронфман), но не все учителя были удачны, а случаи размывания граней тоже постепенно появлялись. Я должен подчеркнуть, что <<Система>> ни до какой степени
не отменяла и не ослабляла традиционные формы обучения.

<<Система>> не вела никакой <<подготовки к вступительным экзаменам>>. Точка зрения была примерно такая: <<Мы учим Природе Вещей, а к этой хренотени отношения не имеем>> (кстати, в рамках <<Системы>> такая <<подготовка>> была технологически невозможна, даже если бы этого очень захотели). 		

У биологов есть свои образцы, и есть возможность подумать над сравнением образовательных технологий. Кстати,  биологический класс Галины Анатольевны Соколовой в 57 школе был организован в 1973~г. по инициативе Константинова.  

\medskip

Стоит напомнить, что обсуждаемые годы были временем расцвета советской математики и временем, когда математическое сообщество предпринимало всевозможные усилия для развития математического образования
в самых разных направлениях. Среди этих действий были успешные (как школы-интернаты в Москве, Ленинграде, Новосибирске и Киеве, основанные в начале 60х), но успехи были далеко не всеобщими - тщательно продуманные проекты, казалось бы обреченные на успех, могли  приводить к результатам, весьма далеким
от желаемых (см., например, \cite{Ner}, \cite{Roh}). Для Константинова было как раз характерно умение находить технологические <<ходы>>%
\footnote{Свои взгляды на преподавание в мат.классах он излагал в своей статье \cite{7}
	1966г. о 7 школе (подписано в печать 24.01.1966). Их статьи ясно, что технология обкатывалась в школе уже по крайней мере третий год, а до того испытывалась в кружках.
	Присоединение 57 школы к <<Cистеме>> - 1968~г., 91 школы - 1969~г., 179 школы -1970 или 71~г. Дату потери 7 школы я не нашел.}%
, находившиеся за пределом круга казавшихся само собой разумеющихся добропорядочных действий. При этом <<Система>> во многом оказывалась и <<школой>>,  и развлечением для работавшей в ней молодежи.

\medskip 

Константинов (все время или не все время - не знаю) вел какой-нибудь мат.класс, а что касается <<Системы>>, то он, по-видимому,
предпочитал косвенное воздействие. Он был вездесущ; приходя  в место, где все работают,  сам садился работать, а также охотно вел оживленные разговоры на произвольные темы (не лишенные, однако, делового содержания). У него были и рычаги влияния, и авторитет, но собиранием власти, как таковой, он, по-моему, не занимался%
\footnote{Он любил цитировать следующую сентенцию из <<Дао Дэ Цзин>>: <<Лучший правитель тот, о котором народ знает лишь то, что он существует.
	Несколько хуже те правители, которые требуют от народа их любить и
	возвышать. Еще хуже те правители, которых народ боится, и хуже всех те,
	которых народ презирает.>> По-моему, он предпочитал, чтобы все шло само собой. Как писали авторы \cite{X}, он <<вовсе не стремился к контролирующей власти или собственной монументальности>>.}... На заднем плане оставались усилия по контактам с внешним миром - каким-то образом все приводилось в состояние <<пахать подано>>. 

Его можно было бы назвать <<энтузиастом>>, но сложившемуся у нас образу <<энтузиаста>> он в одном важном отношении не соответствовал. НикНик хорошо понимал, где грани возможного, и понимал, что область того, что следует желать в такой условно массовой структуре, значительно уже области возможного...   

У участников тусовок обычна склонность превращать тусовку в высшую ценность, такая составляющая присутствовала, и вызывала соответствующую реакцию извне. Но сам Николай Николаевич был настроен именно на обучение школьников математике,
был он человеком абсолютно демократическим, и его дух тогда отражался на образовательной деятельности
тусовки. В школьные кружки могли приходить любые люди, с ними охотно занимались, лишь бы  им была интересна математика. В мат.школы был прием по собеседованию, которое проходило в 4-5 этапов, но человека, которого никто никогда не видел, могли взять сразу, лишь бы он хорошо решал задачи (все прочее от лукавого; к сожалению в более поздние времена, о которых у меня не пойдет речь, многое изменилось). Открытой была и работающая тусовка, лишь бы человек был пригоден, и ему было бы интересно работать.   

\sm

В какой-то момент (кажется, это был 1977~г.) деятели с Мехмата решили отстранить Константинова от руководства Московской математической олимпиады. По причине неформальности связей в тусовке
это было менее очевидной задачей, чем могло показаться затеявшим это лицам. А НикНик,
приговаривая, что <<нельзя всюду навести порядок, как нельзя все заасфальтировать>>,  организовал два новых мероприятия, Турнир Ломоносова (кстати, по  многим разным предметам) и Турнир городов. Они существуют и поныне в виде олимпиад высшего уровня, и, кажется, во многом сохраняют старый олимпиадный дух (впрочем, Конст продолжал ими заниматься). 

Вообще он был очень гибок, без конца изобретал новые формы работы -  не все имело продолжение, но многое получалось; при возникновении проблем переносил деятельность с одного места на другое (благо, что работающая тусовка ни к каким местам привязана не была). Немало он занимался и устройством контактов участников тусовки с людьми из иных (нематематических) миров.

\sm

Осенью 1991~г. по инициативе Николая Николаевича был организован Независимый университет, он же  на начальной фазе его существования обеспечил для оппозиционной мехматской профессуры положение <<Пахать подано>>. Но обучение студентов было уже не его сферой деятельности, и этого с самого начала не предполагалось.

\begin{quotation}
С.И.Комаров: \it ... позвонил Николай Николаевичу. Он, как всегда, без лишних предисловий пригласил меня на следующий
день встретиться... в Моссовете. Удивляюсь. Прихожу. После всех формальностей
с пропусками усаживаемся в большом кабинете за столом под сукном. Появляется
незнакомый мне человек, как оказалось, будущий префект Центрального округа
г.~Москвы А.И.Музыкантский. Смотрим какие-то документы, а Константинов
говорит: <<Вот, это Стаc Комаров, он будет отслеживать все бумаги>>. Это были документы Учредительного Комитета Московского математического университета
(позже он стал Независимым Московским университетом — НМУ). В списке членов Учредительного Комитета на первом месте я увидел фамилию Константинова
с пометкой <<председатель>>, а где-то в середине списка — свою фамилию с пометкой
<<заместитель председателя>>.  Вот так я на следующие пять лет стал менеджером
НМУ.
\end{quotation}

\medskip

Примерно тогда же, по естественной и понятной причине - приходу 90-ых, - прекратила свое существование <<Константиновская система>>. Ее технологические наработки в дальнейшем использовались (и используются), вопрос об этом наследии многогранен и не очевиден.

Так или иначе, Конст продолжил свою разнообразную деятельность, проявляя
положенную ему по штату изобретательность. Я мало видел его после 1992-93~гг., и не берусь рассказывать
об этом периоде его работы.

\begin{center}
	\bf
	\dots\dots\dots
\end{center}

Вернемся к Константинову и ББС. Очевидно, что ему очень понравилось место%
\footnote{Там очень красивые и разнообразные окрестности.}, ему понравился Директор Перцов,
и, главное, он счел (и был в этом глубоко прав), что тамошние работы, тамошние прогулки
 и тамошнее общение  полезны для общего развития математически ориентированных молодых людей.  
И вот с 1969~г. начали появляться стройотряды, в которых основную часть составляли
<<математики>> - ученики и выпускники физико-математических школ, а также  студенты и аспиранты математических специальностей (эта социальная группа  преобладала в мужской половине отрядов, руководство  тоже чаще состояло из ее представителей). Обычным временем заезда таких отрядов были июль, август, Новый год (короткий промежуток между зачетной и экзаменационной сессией) и зимние студенческие каникулы (были еще <<физтехи>>, которые представляли свою спаянную компанию, они ездили, кажется, в марте-апреле). В опубликованных
списках
стройотрядов (в них могут быть лакуны) Константинов упоминается в  
июле 1967~г.,
июле 1969~г.,
июле 1972~г.,
июле 1973~г.,
июле-августе 1977~г.,
зимой 1979~г.,
июне-июле 1979~г. Вроде бы это не очень много (да и появлялся он не на полные сроки), но  надо иметь в виду, что на биостанцию приезжали люди, связанные между
собой <<в другой жизни>>. И связанные между собой не только знакомствами, но и совместной работой в <<Системе>>.
Фактически, и сами эти стройотряды стали своего рода частью <<Системы>>.

Люди приезжали за разным - тусоваться, отдыхать, за природой, за романтикой... (напомню, работали
беззарплатно, обеспечивались харчи и оплачивались железнодорожные билеты; выдавалась рабочая одежда). Летом был восьмичасовой рабочий день (без дураков), за исключением пилорамы, где считалось, что 8 часов должны работать основные станки. Дальше надо было переработать обрезки на дрова, оприходовать согласно нормам складирования всю образовавшуюся продукцию, убрать помещение, подготовить станки. Что занимало никак не меньше еще двух часов...  Работы были разные, 
были хорошие (например, плотницкие или пилорамные - на пилораму  люди рвались, несмотря на ненормированный рабочий день), а были, скажем, небезынтересные земляные работы - целинную зеленоватую глину, пропитанную камнями любых размеров, не брал никакой ручной инструмент, кроме кайла. Было много разного <<труда, освобожденного от разума>>. Кстати, Конст в свой первый приезд начал с  разделения щебня на фракции разного размера на устройстве с характерным названием <<грохот>>. Была тьма внезапно возникавших работ в связи с  естественными или искусственными кризисами (например, промерзание водопроводного ручья или наезд санэпидемстанции) или в связи с очередной поломкой, которые иногда шли бесконечной чередой  (помнится, на собрании в Москве, где в частности обсуждались предстоявшие работы на озере Бульдозере%
\footnote{Котлован под водоем, названный так в 1978г. в честь застрявшего там бульдозера.},
Перцов сказал, что было бы хорошо купить вертолет, на что сидевший сзади Конст заметил: <<Представляете,
вертолет застрял в облаках>>). Вообще о тамошних работах см. много романтических и не очень романтических воспоминаний в \cite{Kal}.

Люди были разные, работали по-разному, но, в общем,  работали. В стройотряд могли попасть далеко не все желающие, а для школьников поездка на ББС была чем-то вроде высокой награды от учителей, положительный отбор в отряд был. Ну и главное, что для части тусовки строительство ББС стало таким же общим делом (Res publica) как, скажем, ведение классов или организация олимпиад (спектр участников был несколько иным, но в <<Системе>> однородность и не предполагалась). Злые языки, конечно, вспомнят  присутствовавшую в стройотряде идею, что биостанция нужна для того, чтобы ее строить. Но и сам ее директор был не чужд таких настроений.

\begin{figure}
	\includegraphics[width=157mm]{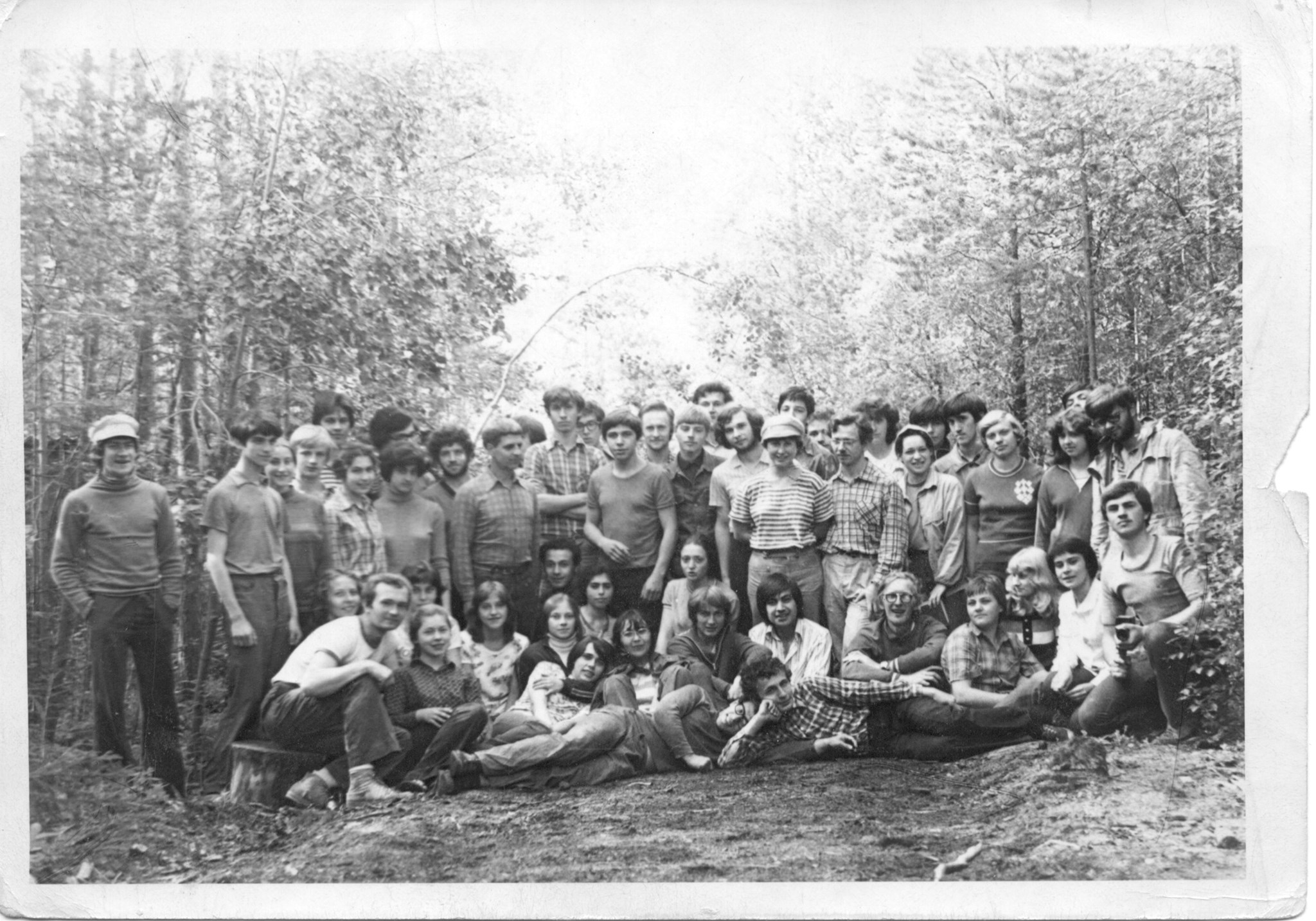}
	\caption{ББС МГУ. Стройотряд. Июль, 1979~г. Фото С.Л.Мехедова. Конст - под нависшей березой без листьев. Видна степень многолюдности типового летнего стройотряда.}
	
	 Укажем людей, впоследствии ставших известными профессионалами (по названиям профессий открываются ссылки). Крайний слева - С.Ю.Оревков 
		(\href{https://mathscinet.ams.org/mathscinet/MRAuthorID/202757   target="_blank"}{математик}),
		сидящий крайний справа - 
		Д.С.Богданов 
		(\href{http://www.isa.ru/index.php?option=com_content&view=article&id=301&Itemid=122&lang=ru   target="_blank"}{программист} и 
		\href{http://www.bards.ru/person.php?id=2672  target="_blank"} {гитарист}), рядом с Констом справа - 
		М.М.Капранов 
		(\href{https://mathscinet.ams.org/mathscinet/MRAuthorID/200368   target="_blank"}{математик}), лежит справа А.М.Левин 
		(\href{https://mathscinet.ams.org/mathscinet/MRAuthorID/291654   target="_blank"}{математик}), сидит над его головой 
		Н.Х.Агаханов 
		(\href{http://www.mathnet.ru/php/person.phtml?&personid=37358%20target=_blank&option_lang=eng}{деятель школьных математических олимпиад}), 
		стоит четвертый слева - А.Л.Городенцев 
		(\href{https://mathscinet.ams.org/mathscinet/MRAuthorID/241278   target="_blank"}{математик}),
		стоит выше девушки в тельняшке 
		Ан.Резников 
		(\href{https://mathscinet.ams.org/mathscinet/MRAuthorID/333309   target="_blank"}{математик}). 
		За спиной М.М.Капранова справа - 
		А.А.Васильев 
		(\href{http://www.iki.rssi.ru/books/2012vasiliev_ar.pdf}{физик}),
		много ездивший работать на биостанцию в смутные годы, слева от Конста -  
		В.В.Прасолов 
		(\href{https://mathscinet.ams.org/mathscinet/MRAuthorID/210806   target="_blank"}{автор книг по математике}),
		третий стоящий слева от девушки в тельняшке - 
		П.М.Ахметьев(\href{https://mathscinet.ams.org/mathscinet/MRAuthorID/226682   target="_blank"}{математик}),
		сидит слева А.П.Романов. Сидит под М.М.Капрановым - 
		П.Я.Грозман 
		(\href{https://mathscinet.ams.org/mathscinet/MRAuthorID/609811  target="_blank}{математик}), лежит слева 
		В.Н.Шандер 
		(\href{https://mathscinet.ams.org/mathscinet/MRAuthorID/192823   target="_blank"}{математик}).
		Человек в клечатой рубахе (около девушки в тельняшке) - 
		Ю.П.Лысов (один из важнейших соратников Конста, впоследствии 
		\href{https://scholar.google.com/scholar?hl=ru&as_sdt=0%2C5&q=lysov+y&btnG=  target="_blank"}{молекулярный биолог}).
		Слева от девушки в тельняшке - В.О.Бугаенко
		(\href{http://www.mathnet.ru/php/person.phtml?&personid=19977%20   target=_blank&option_lang=eng}
			{деятель школьного образования}, много сотрудничал с Констом).
		Над ногами лежащего Левина - 
		В.В.Серганова 
		(\href{https://mathscinet.ams.org/mathscinet/MRAuthorID/158860   target="_blank}{математик}). Второй слева - 
		М.З.Шапиро 
		(\href{https://mathscinet.ams.org/mathscinet/MRAuthorID/249594   target="_blank"}{математик}).
		В списке стройотряда \cite{Kal} присутствуют также 
		В.А.Васильев 
		(\href{https://mathscinet.ams.org/mathscinet/MRAuthorID/327216   target="_blank"}{математик}),
		А.Ю.Вайнтроб 
		(\href{https://mathscinet.ams.org/mathscinet/MRAuthorID/204643   target="_blank}{математик}),
		С.А.Кулешов(\href{https://mathscinet.ams.org/mathscinet/MRAuthorID/249337    target="_blank"}{математик}),
		Е.И.Полетаева (\href{https://mathscinet.ams.org/mathscinet/MRAuthorID/316083    target="_blank"}{математик}),
		но я что-то не могу найти их на фото. Фотограф -- 
		\href{https://scholar.google.com/scholar?hl=ru&as_sdt=0%2C5&q=mekhedov+sergei&oq=mekhedov9  target="_blank"}
		{молекулярный биолог}.
		Я, скорее всего, в этот момент тоже фотографировал (но менее удачно). Состав отряда не вполне типичен, обычно там было больше представителей (представительниц) биологического мира.

\end{figure}

Стоит заметить, что работа стройотряда требовала больших  затрат времени и сил директора,
Николая Андреевича Перцова. Он, владея на профессиональном уровне  несколькими строительными специальностями, должен был много что объяснять все новым и новым людям, хоть и настроенным на добросовестность, но ничего не умеющим. С КПД стройотряда были проблемы, но постройки он за собой оставлял, а как показали итоги десятилетия, когда все было в состоянии свободного развала (90-ые), биостанция строилась хорошо%
\footnote{Скажем несколько слов о дальнейшем (для математиков, строивших биостанцию в 60-80х годах, если им этот текст попадется на глаза). Хотя на биостанции очень многое было завязано на личность Перцова, следующие несколько лет благодаря нескольким дельным сотрудникам там, видимо, продолжалась прежняя жизнь. В частности, от этого времени остались новые постройки (хотя Перцов бы их не одобрил).
	\newline
	А потом пришли 90е.  За долги было отключено электричество, дальше группа местных жителей зимой спилила столбы ЛЭП и <<собрала урожай>>  - провода - для сдачи их на лом из цветного металла  (классический для 90х мелкий бизнес). Менее <<ценная>> телефонная линия уцелела (впрочем, спиливание столбов при еще работающем телефоне могло бы быть пресечено). Технический персонал биостанции, получая финансирование от Московского университета, начал действовать и осваиваивать средства и возможности биостанции в своих интересах.		
	\newline
	Биостанция была превращена в Перцовым в место, нужное для биофака и многих биологов, что сильно способствовало ее выживанию. Иногда на станцию заезжали представители старого стройотряда, которые проводили там <<неотложные аварийно-восстановительные работы>>. Основную роль играл А.В.Андрианов (физик с физфака МГУ), 
	имевший местную кличку <<мозг>>, который, в частности, ежегодно приводил в порядок слегка питавший станцию электричеством дизель. Что касается  станционной инфраструктуры, то она в условиях развала, воровства и бардака оказалось очень живучей и дотянула до 2005г.
	\newline 
	Быстрому ренессансу способствовало обстоятельство, отчасти, случайное. В 2002г. на деньги, выбитые одним известным биологом, была куплена новая посудина, капитаном на ней стал живший в Пояконде на пенсии  Л.~Д.~Папунашвили, которого тогда представляли как  <<бывшего капитана дальнего плавания>> (в реальности капитан первого ранга, командир атомной подлодки и прочая). <<Давыдыч>> быстро привел в порядок станционные суда и, видимо, никого не трогал и его никто не трогал. В 2005г. очередному заведующему базой  хватило разумения на него <<наехать>>.  Спустя короткое время Давыдыч стал заведующим базой (хозяйственная должность), и быстро вслед за этим директором (научная должность) стал А.Б.Цетлин. С этими двумя назначениями несколько затянувшиеся в рамках одной отдельно взятой биостанции <<святые девяностые>> были выключены  как поворотом тумблера.}. Со смертью Перцова и концом <<Системы>> существование стройотрядов в прежнем виде стало невозможным.

\bigskip

Так или иначе,
в начале 70-х 
 в стройотряде на ББС вскоре все пошло своим чередом (что соответствовало  нравам тусовки), а Константинов затеял новый проект. В 1974~г., договорившись
с начальством какого-то колхоза   в  Эстонии%
\footnote{Кстати (см. \cite{Gor}) в том же 1974~г. Конст основал известную и работающую поныне школьную билогическую базу в Ковде (поморское село в Мурманской области).
Н.~А.~Горяшко (см. \cite{Gor}):	
{\it Летом 1974 г. группа биологов-девятиклассников во главе с Н.~Н.~Константиновым выехала на Белое море, чтобы подыскать подходящее место для постоянной практики биокласса.  << Сама я }[цитируется Г.~А.~Соколова]{\it  в тот год поехать не смогла, а поехал Николай Николаевич Константинов. Они сначала приехали в Пояконду, там купили лодку. Потом они на лодке плыли вдоль берега Белого моря, и нашли эту деревню, здесь нашли этот дом. Заплатили хозяевам 10 рублей, хозяева были этим глубоко потрясены. Раньше у них просто сбивали замок, заходили и жили. Я приехала на следующий год, мы опять заплатили 10 рублей, и у нас был дом, и была лодка... В 76 г. мы этот дом купили>>}.%
},   он
устроил там летний математический лагерь%
\footnote{См. рассказ посетившей лагерь учительницы литературы  в \cite{Kle}. 
	{\it <<В летнем математическом лагере по вечерам у костра проводились беседы обо всем на свете, от индийской философии до проблем генетики>>...<<Никаких слов о том, «что такое хорошо», Константинов и его помощники не произносят.>>}
}. Мероприятие (<<Констлагерь>>) было в несколько ином стиле - там была математическая программа и туда созывались
школьники и студенты из самых разных городов страны. Там тоже велись определенные работы, как
для обеспечения жизни лагеря, так и для того, чтобы у личного состава было поменьше свободного времени.
Помню, например, сцену, когда 17 человек на строительстве объекта <<Чайхана>> в процессе бурных обсуждений и трудовых подвигов совместными усилиями прибили за  рабочий день одну доску. Лишь с некоторым запаздыванием из нескольких реплик я понял, что Конст  развлекался, управляя этим процессом, и считал, что в данных
условиях эта версия  <<бега на месте>> полезна, а нормально понимаемая работа могла бы стать
(и становилась) источником сложностей.

Лагерь  просуществовал до 1989г., на следующее лето Эстония Конста не впустила.

\begin{center}
	\bf
	\dots\dots\dots
\end{center}

\begin{figure}
	\includegraphics[width=70mm]{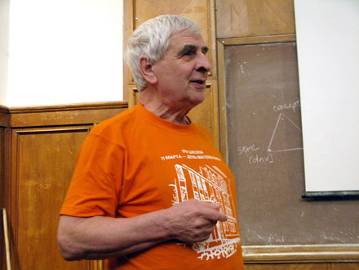}
	\caption{
		Конст у доски. По виду, это одна из поточных аудиторий в корпусах МГУ, построенных в 1953г.
		Он часто читал  публичные лекции для школьников и с легкостью работал в шумящих аудиториях,
		раз за разом включая шум или буйство и изящно выключая их, как только это становилось необходимым. Когда-то по окончании Московской олимпиады он выступал в переполненных аудиториях
		01 или 02 в Шайбе (самых больших аудиториях тогдашнего МГУ). Степень мотивированности собравшихся была большой, но не полной - аудитория не была безобидной - и зрелище применяемой НН техники  контроля  было забавным и поучительным.}
\end{figure}

НН учился на физфаке МГУ. Защитил диссертацию по математике 
{\it Некоторые задачи теоретико-множественной геометрии плоских кривых}, МГУ, 1963.
База MathSciNet показывает у него 2 математических работы с одним цитированием (что
современные научные джентльмены, по-видимому, должны считать предосудительным). 

Еще у него была (невидимая никакими скопусами) пионерская работа по прикладной математике, см. \cite{Kosh}--\cite{Lev2}.
В 1968~г., работая  в лаборатории А.С.Кронрода в ИТЭФ, он  со товарищи
сделал фильм <<{\it  Кошечка}>> \cite{Kosh}.
Это был первый опыт по компьютерной мультипликации в СССР и один из первых в мире. Причем <<Кошечка>>
надолго осталась единственной такой работой, в которой движение вычислялось с помощью дифференциальных уравнений (более простой для мультипликации путь - интерполяция, художник рисует несколько положений объекта, а машина вычисляет
промежуточные состояния). Дифференциальные уравнения должны были выдавать движения, которые визуально похожи на настоящие,
задача описывать реальное движение животного  не ставилась. Коты и кошки наотрез отказывались 
позировать для решения научной задачи, поэтому (как рассказывал НН) одному из участников проекта -- \href{https://rg.ru/2003/11/19/vystavka.html}{В.В.Минахину} -- пришлось много ползать на четвереньках, размышляя над описанием движений тех или иных участков конечностей. Кстати, работа делалась на ламповых ЭВМ (БЭСМ-4), каковые были несколько менее удобны, чем современные компы.

К сожалению, <<Кошечка>> не повлекла практических приложений, а исследования в этом направлении не продолжались. Не привлекли тогда подобные задачи и биологов. Но
широко известная в узких кругах <<Кошечка>>, по-видимому, имела косвенное воздействие, способствовав появлению  у нескольких талантливых молодых математиков  интереса к математическим вопросам биологии.

\begin{center}
	\bf
	\dots\dots\dots
\end{center}

 Значение Константинова для советской/российской математики последних 55 лет невозможно
ни оценить, ни описать - без него окружающий математиков мир просто был бы совсем иным. 
За  тридцатилетие 1960-1990~гг., о котором я писал, деятельность Константинова и <<Системы>> коснулась 
многих десятков тысяч человек (и совсем не только будущих профессиональных математиков)... Очень  трудно оценить значение дальнейших <<кругов на воде>> и
значение его побочных
инициатив, которые, бывало, претворялись в жизнь  (как создание биологических классов). 

Был он человеком абсолютного благородства, и был он человеком удивительно мудрым. В своей работе
ему было суждено избежать больших трагических ошибок, столь обычных на пути крупных общественных
деятелей, имевших самые лучшие намерения. И, наконец, встреча лично с ним была поворотной точкой в судьбе многих талантливых молодых людей...

\noindent
\tt University of Vienna, Department of Mathematics;
\\
Институт теоретической и экспериментальной физики (ИТЭФ);
\\
Мехмат МГУ;
\\
Институт проблем передачи информации РАН.
\\
e-mail:yurii.neretin@univie.ac.at
\\
URL:mat.univie.ac.at/$\sim$neretin

 \begin{figure}
	\includegraphics[width=157mm]{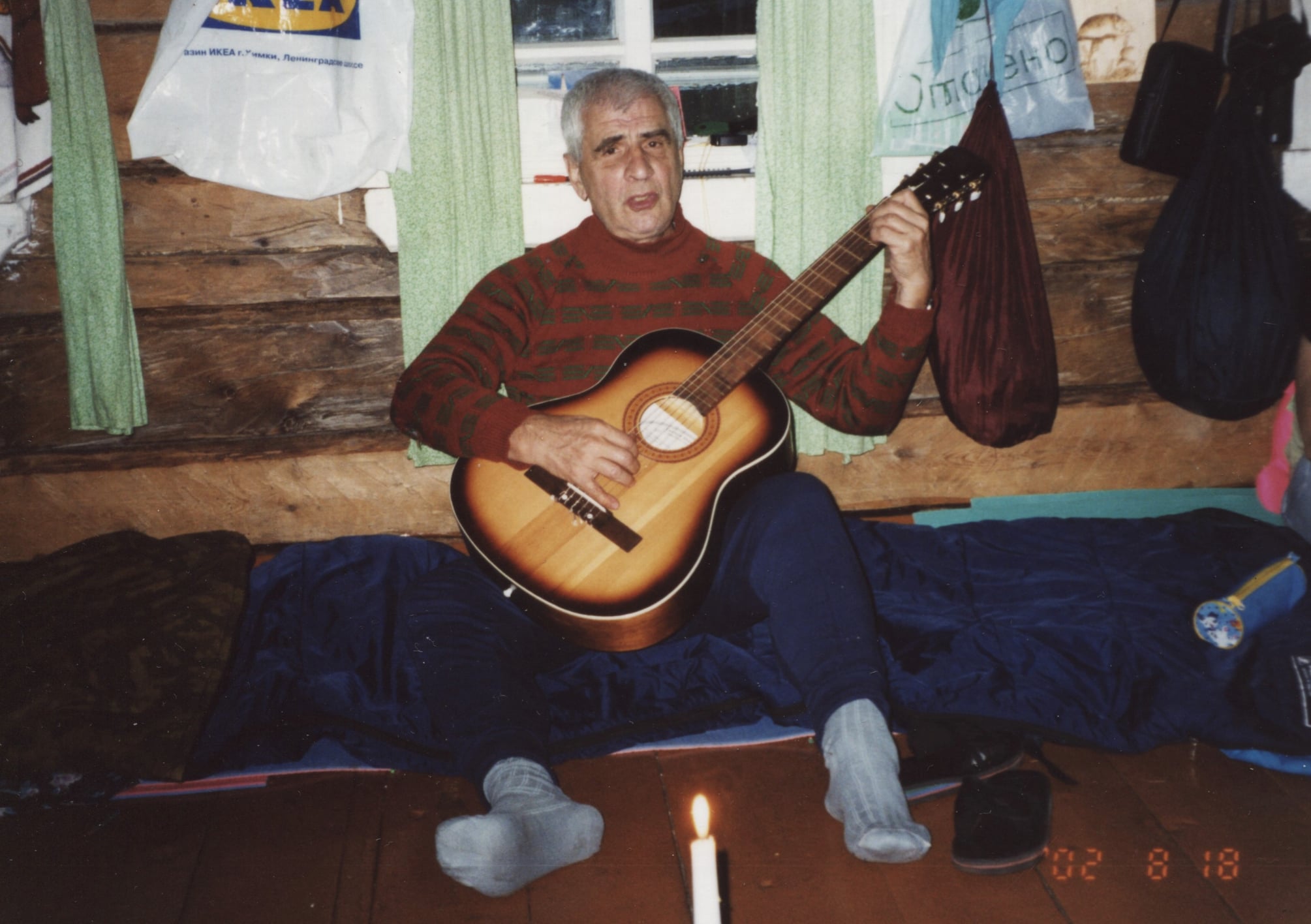}
	\caption{Школьная биологическая база в Ковде, август 2002}
\end{figure}

\end{document}